\documentclass[a4paper,fleqn]{article}
\usepackage[english]{babel}
\usepackage{amsmath,amssymb,amsfonts,amsthm}
\usepackage{mathrsfs}
\usepackage{graphicx}
\usepackage[numbers,sort&compress]{natbib}

\newtheorem{theorem}{Theorem}
\newtheorem{lemma}{Lemma}
\newtheorem{example}{Example}
\newtheorem*{barthm}{Barabanov's Theorem}
\newtheorem*{wirthm}{Wirth's Theorem}

\newcommand{\ecc}{\mathop{\mathrm{ecc}}}
\newcommand{\conv}{\mathop{\mathrm{conv}}}
\newcommand{\setA}{\mathscr{A}}
\newcommand{\setB}{\mathscr{B}}
\newcommand{\setP}{\mathscr{P}}
\newcommand{\setK}{\mathscr{K}}
\newcommand{\setI}{\mathscr{I}}

\DeclareMathOperator*{\infx}{\vphantom{sup}inf}

\begin{document}
\date{}

\title{An explicit Lipschitz constant\\ for the
joint spectral radius\thanks{Supported by the Federal Agency
for Science and Innovations of Russian Federation, contract no.
02.740.11.5048, and by Russian Foundation for Basic Research,
project no. 10-01-00175.}}

\author{Victor Kozyakin\\[5mm]
Institute for Information Transmission Problems\\ Russian
Academy of Sciences\\ Bolshoj Karetny lane 19, Moscow 127994
GSP-4, Russia}

\maketitle

\begin{abstract}
In 2002 F.~Wirth has proved that the joint spectral radius of
irreducible compact sets of matrices is locally Lipschitz
continuous as a function of the matrix set. In the paper, an
explicit formula for the related Lipschitz constant is
obtained.

\medskip

\medskip\noindent PACS number 02.10.Ud; 02.10.Yn


\medskip

\end{abstract}

\section{Introduction}\label{S-intro}
Information about the rate of growth of matrix products with
factors taken from some matrix set is of great importance in
various problems of control theory
\cite{BrayTong:TCS80,Bar:AIT88-2:e,BerWang:LAA92}, wavelet
theory \cite{DaubLag:SIAMMAN92,DaubLag:LAA92,ColHeil:IEEETIT92}
and other fields of mathematics. One of the most prominent
values characterizing the exponential rate of growth of matrix
products is the so-called joint or generalized spectral radius.

Let $\mathbb{K}=\mathbb{R},\mathbb{C}$ be the field of real or
complex numbers, and $\setA\subset \mathbb{K}^{d\times d}$ be a
set of $d\times d$ matrices. As usual, for $n\ge1$ denote by
$\setA^{n}$ the set of all $n$-products of matrices from
$\setA$; $\setA^{0}=I$.

Given a norm $\|\cdot\|$ in ${\mathbb{K}}^{d}$, the limit
\begin{equation}\label{E-JSRad}
\rho({\setA})=
\limsup_{n\to\infty}\|\setA^{n}\|^{1/n},
\end{equation}
with $\|\setA^{n}\|=\sup_{A\in\setA^{n}}\|A\|$ is called
\emph{the joint spectral radius} of the matrix set $\setA$
\cite{RotaStr:IM60}. The limit in (\ref{E-JSRad}) is finite for
bounded matrix sets $\setA\subset \mathbb{K}^{d\times d}$ and
does not depend on the norm $\|\cdot\|$. As shown in
\cite{RotaStr:IM60}, for any $n\ge 1$ the estimates
$\rho({\setA})\le \|\setA^{n}\|^{1/n}$ hold, and therefore the
joint spectral radius can be defined also by the following
formula:
\begin{equation}\label{E-JSR}
\rho({\setA})=
\inf_{n\ge 1}\|\setA^{n}\|^{1/n}.
\end{equation}
If $\setA$ is a singleton set then (\ref{E-JSRad}) turns into
the known Gelfand formula for the spectral radius of a linear
operator. By this reason sometimes (\ref{E-JSRad}) is called
the generalized Gelfand formula \cite{ShihWP:LAA97}.

Besides (\ref{E-JSRad}) and (\ref{E-JSR}), there are quite a
number of different equivalent definitions of $\rho({\setA})$
in which the norm in (\ref{E-JSRad}) is replaced by the
spectral radius
\cite{DaubLag:SIAMMAN92,DaubLag:LAA92,Els:LAA95} or the trace
of a matrix \cite{ChenZhou:LAA00}, or by a uniform non-negative
polynomial of even degree \cite{ParJdb:LAA08}. Sometimes
$\rho({\setA})$ is defined in terms of existence of specific
norms \cite{Bar:AIT88-2:e,Prot:FPM96:e} (the Barabanov and
Protasov norms). Unfortunately, the common feature of all the
mentioned definitions is their nonconstructivity. In all these
definitions the value of $\rho(\setA)$ is defined either as a
certain limit or as a result of some ``existence theorems'',
which essentially complicates the analysis of properties of the
joint spectral radius.

In the paper, we are concerned with properties of the joint
spectral radius $\rho({\setA})$ as a function of $\setA$ for
compact (i.e. closed and bounded) matrix sets $\setA$. In this
case it is convenient to denote the set of all nonempty bounded
subsets of $\mathbb{K}^{d\times d}$ by
$\setB(\mathbb{K}^{d\times d})$, and the set of all nonempty
compact subsets of $\mathbb{K}^{d\times d}$ by
$\setK(\mathbb{K}^{d\times d})$. Both of these sets become
metric spaces if to endow them with the usual Hausdorff metric
\[
H(\setA,\setB):=
\max\left\{\sup_{A\in\setA}\infx_{B\in\setB}\|A-B\|,
~\sup_{B\in\setB}\infx_{A\in\setA}\|A-B\|\right\}.
\]
In doing so the space $\setK(\mathbb{K}^{d\times d})$ is proved
to be complete while the set $\setI(\mathbb{K}^{d\times d})$ of
all irreducible compact matrix families is open and dense in
$\setK(\mathbb{K}^{d\times d})$.

In 2002 F.~Wirth has proved \cite[Cor.~4.2]{Wirth:LAA02} that
the joint spectral radius of irreducible compact matrix sets
satisfies the local Lipschitz condition.
\begin{wirthm}
For any compact set $\setP\subset\setI(\mathbb{K}^{d\times d})$
there is a constant $C$ (depending on $\setP$ and the norm
$\|\cdot\|$ in $\mathbb{K}^{d\times d}$) such that
\[
\left|\rho(\setA)-\rho(\setB)\right|\le C\cdot H(\setA,\setB),\quad
\forall~ \setA,\setB\in\setP.
\]
\end{wirthm}
The aim of the present paper is to obtain an explicit
expression for the constant $C$ in the above inequality.

As demonstrated the following example the joint spectral radius
is not locally Lipschitz continuous if to discard supposition
about irreducibility of a matrix set.

\begin{example}\label{Ex1}\rm
Consider the matrix set $\setA_{\varepsilon}$ composed of a
single matrix
\[
A_{\varepsilon}=\left(\begin{array}{cc}
1 & 1\\ \varepsilon &
1\end{array}\right),
\]
depending on the real parameter $\varepsilon>0$.

Clearly, the singleton matrix set $\setA_{0}$ is not
irreducible. Besides,
$\rho(A_{\varepsilon})=1+\sqrt{\varepsilon}$, and therefore
\[
\left|\rho(\setA_{\varepsilon})-\rho(\setA_{0})\right|=
\left|\rho(A_{\varepsilon})-\rho(A_{0})\right|=\sqrt{\varepsilon},
\]
whereas
$H(\setA_{\varepsilon},\setA_{0})=\|A_{\varepsilon}-A_{0}\|=\varepsilon
c$, where $c$ is some constant.
\end{example}

\section{Main result}\label{S-irred}

Given a matrix set $\setA\subset \mathbb{K}^{d\times d}$, for
each $n\ge 1$ denote by ${\setA}_{n}$ the set of all
$k$-products of matrices from $\setA\bigcup \{I\}$ with $k\le
n$, that is ${\setA}_{n}=\cup_{k=0}^{n}\setA^{k}$. Denote also
by ${\setA}_{n}(x)$ the set of all the vectors $Ax$ with
$A\in{\setA}_{n}$. Let $\|\cdot\|$ be a norm in
$\mathbb{K}^{d}$ then $\mathbb{S}(t)$ stands for the ball of
radius $t$ in this norm.

Let us call the \emph{$p$-measure of irreducibility} of the
matrix set $\setA$ (with respect to the norm $\|\cdot\|$) the
quantity $\chi_{p}({\setA})$ determined as
\[
\chi_{p}({\setA}) = \inf_{\substack{x\in\mathbb{R}^{d}\\ \|x\|=1}}
\sup \{ t:\mathbb{S}(t)\subseteq{\conv} \{ {\setA}_{p}(x)\cup {\setA}_{p}(-x)\}
\}.
\]

Under the name `the measure of quasi-controllability' the
measure of irreducibility $\chi_{p}({\setA})$  was introduced
and investigated in
\cite{KozPok:DAN92:e,KozPok:CADSEM96-005,KKP:CESA98} where the
overshooting effects for the transient regimes of linear remote
control systems were studied. The reason why the quantity
$\chi_{p}({\setA})$ got the name `the measure of
irreducibility' is in the following lemma.

\begin{lemma}\label{L-qcontr}
Let $p\ge d-1$. The matrix set $\setA$ is irreducible if and
only if $\chi_{p}({\setA})> 0$.
\end{lemma}

The proof of Lemma~\ref{L-qcontr} can be found in
\cite{KozPok:CADSEM96-005,KKP:CESA98}. In these works it is
proved also that, for compact irreducible matrix sets, the
quantity $\chi_{p}({\setA})$ continuously depends on $\setA$ in
the Hausdorff metric.

\begin{theorem}\label{T-main}
For any pair of matrix sets $\setA\in\setI(\mathbb{K}^{d\times
d})$, $\setB\in\setB(\mathbb{K}^{d\times d})$ for each $p\ge
d-1$ it is valid the inequality
\begin{equation}\label{E-mainineq}
\left|\rho(\setA)-\rho(\setB)\right|\le
\nu_{p}(\setA) H(\setA,\setB),
\end{equation}
where
\[
\nu_{p}({\setA})=
\frac{\max\{1,\|\setA\|^{p}\}}{\chi_{p}({\setA})}.
\]
\end{theorem}

Taking into account that the quantity $\nu_{p}(\setA)$
continuously depends on $\setA$ in the Hausdorff metric, and
hence it is bounded on any compact set
$\setP\subset\setI(\mathbb{K}^{d\times d})$,
Theorem~\ref{T-main} implies Wirth's Theorem. However, unlike
to Wirth's Theorem, in Theorem~\ref{T-main} neither compactness
nor irreducibility of the matrix set $\setB$ is assumed.

As will be shown below under the proof of Theorem~\ref{T-main},
in fact even more accurate estimate than (\ref{E-mainineq})
holds:
\[
\left|\rho(\setA)-\rho(\setB)\right|\le
\frac{\max\{1,(\rho(\setA))^{p}\}}{\chi_{p}({\setA})}H(\setA,\setB).
\]
However, this last estimate is not quite satisfactory because
practical evaluation of the quantity $\rho(\setA)$ is a
problem. At the same time the quantity $\nu_{p}({\setA})$ in
(\ref{E-mainineq}) \textsl{can be evaluated in a finite number
of algebraic operations involving only information about
$\setA$}.

Remark also that whereas the value of the joint spectral radius
is independent of a norm in $\mathbb{K}^{d\times d}$, the
quantities $\nu_{p}(\setA)$, $\chi_{p}(\setA)$ and
$H(\setA,\setB)$ in (\ref{E-mainineq}) do depend on the choice
of the norm $\|\cdot\|$ in $\mathbb{K}^{d\times d}$.

At last, point out that in the case when both of the matrix
sets $\setA$ and $\setB$ are irreducible and compact, that is
$\setA,\setB\in\setI(\mathbb{K}^{d\times d})$, inequality
(\ref{E-mainineq}) can be formally strengthened:
\[
\left|\rho(\setA)-\rho(\setB)\right|\le
\min\left\{\nu_{p}(\setA),\nu_{p}(\setB)\right\} H(\setA,\setB).
\]

\section{Auxiliary statements}\label{S-aux}

To prove Theorem~\ref{T-main} we will need some auxiliary
notions and facts related to the irreducible matrix sets. The
principal technical tool in proving Theorem~\ref{T-main} will
be the notion of the Barabanov norm mentioned above, existence
of which follows from the next theorem~\cite{Bar:AIT88-2:e}.

\begin{barthm}
The quantity $\rho$ is the joint (generalized) spectral radius
of the matrix set $\setA\in\setI(\mathbb{K}^{d\times d})$ if
and only if there is a norm $\|\cdot\|_{b}$ in
${\mathbb{K}}^{d}$ such that
\begin{equation}\label{Eq-mane-bar}
\rho\|x\|_{b}\equiv
\max_{A\in\setA}\|Ax\|_{b}.
\end{equation}
\end{barthm}

In what follows a norm satisfying (\ref{Eq-mane-bar}) is called
\emph{a Barabanov norm} corresponding to the matrix set
$\setA$.

In the next elementary lemma, a simple way to get both upper
and lower estimates for the joint spectral radius is suggested.

\begin{lemma}\label{L-rhorel}
Let $\setA$ be a nonempty matrix set from $\mathbb{K}^{d\times
d}$. If, for some $\alpha$,
\begin{equation}\label{E-alphacond}
\sup_{A\in\setA}\|Ax\|\le\alpha\|x\|,\quad\forall~x\in\mathbb{K}^{d},
\end{equation}
then $\rho(\setA)\le\alpha$. If, for some $\beta$,
\begin{equation}\label{E-betacond}
\sup_{A\in\setA}\|Ax\|\ge\beta\|x\|,\quad\forall~x\in\mathbb{K}^{d},
\end{equation}
then $\rho(\setA)\ge\beta$.
\end{lemma}

\begin{proof}
Clearly, the constants $\alpha$ and $\beta$ may be thought of
as non-negative. To prove the first claim note that
(\ref{E-alphacond}) implies the inequality
$\|\setA\|=\sup_{A\in\setA}\|A\|\le\alpha$. Then
$\|\setA^{n}\|=\sup_{A_{i}\in\setA}\|A_{n}\cdots
A_{2}A_{1}\|\le\alpha^{n}$, and $\rho(\setA)\le\alpha$ by the
definition (\ref{E-JSRad}).

Similarly, to prove the second claim note that
(\ref{E-betacond}) implies, for each $n=1,2,\dotsc$, the
inequality
\begin{multline*}
\sup_{A_{i}\in\setA}\|A_{n}\cdots A_{2}A_{1}x\|
=\sup_{A_{1}\in\setA}\sup_{A_{2}\in\setA}\ldots\sup_{A_{n}\in
\setA}\|A_{n}\cdots A_{2}A_{1}x\|\\
\ge\beta^{n}\|x\|,\quad\forall~x\in\mathbb{K}^{d}.
\end{multline*}
Hence $\sup_{A_{i}\in\setA}\|A_{n}\cdots
A_{2}A_{1}\|\ge\beta^{n}$. Then
$\|\setA^{n}\|=\sup_{A_{i}\in\setA}\|A_{n}\cdots
A_{2}A_{1}\|\ge\beta^{n}$, and $\rho(\setA)\ge\beta$ by the
definition (\ref{E-JSRad}). The lemma is proved.
\end{proof}

Following to \cite{Wirth:CDC05}, for convenience of comparison
of norms in ${\mathbb{K}}^{d}$ let us introduce an appropriate
notion. Given two norms $\|\cdot\|'$ and $\|\cdot\|''$ in
${\mathbb{K}}^{d}$, define the quantities
\begin{equation}\label{E-eccentr}
e^{-}(\|\cdot\|',\|\cdot\|'')=\min_{x\neq0}\frac{\|x\|'}{\|x\|''},\quad
e^{+}(\|\cdot\|',\|\cdot\|'')=\max_{x\neq0}\frac{\|x\|'}{\|x\|''}.
\end{equation}

Since all norms in ${\mathbb{K}}^{d}$ are equivalent then the
quantities $e^{-}(\cdot)$ and $e^{+}(\cdot)$ are well defined,
and
\[
0< e^{-}(\|\cdot\|',\|\cdot\|'')\le
e^{+}(\|\cdot\|',\|\cdot\|'')< \infty.
\]
Therefore the quantity
\begin{equation}\label{E-defeccentr}
\ecc(\|\cdot\|',\|\cdot\|'')=
\frac{e^{+}(\|\cdot\|',\|\cdot\|'')}{e^{-}(\|\cdot\|',\|\cdot\|'')}\ge 1,
\end{equation}
called \textit{the eccentricity} of the norm $\|\cdot\|'$ with
respect to the norm $\|\cdot\|''$, is also well defined.

\section{Proof of Theorem~\ref{T-main}}\label{S-main}

We will prove Theorem~\ref{T-main} in two steps. First,
slightly modifying the idea of the proof from
\cite[Cor.~4.2]{Wirth:LAA02}, we will show in
Section~\ref{S-step1} that under the conditions of
Theorem~\ref{T-main} the eccentricity of any Barabanov norm
$\|\cdot\|_{\setA}$ for the matrix set $\setA$ with respect to
the norm $\|\cdot\|$ may serve as the Lipschitz constant for
the joint spectral radius, that is
\begin{equation}\label{E-Lipecc}
\left|\rho(\setA)-\rho(\setB)\right|\le
\ecc(\|\cdot\|_{\setA},\|\cdot\|) H(\setA,\setB).
\end{equation}
Then, using the techniques of the measures of irreducibility
(see, e.g., \cite{KozPok:DAN92:e,KKP:CESA98,Koz:ArXiv08-2}), we
will prove in Section~\ref{S-step2} the estimate
\begin{equation}\label{E-nubound}
\ecc(\|\cdot\|_{\setA},\|\cdot\|)\le \nu_{p}({\setA}):=
\frac{\max\{1,\|\setA\|^{p}\}}{\chi_{p}({\setA})}.
\end{equation}

\subsection{Proof of estimate (\ref{E-Lipecc})}\label{S-step1}

Let $\|\cdot\|_{\setA}$ be some Barabanov norm for the matrix
set $\setA$. By definition of the Hausdorff metric, for any
matrix $B\in\setB$ there is a matrix $A_{B}\in\setA$ such that
$\|B-A_{B}\|\le H(\setA,\setB)$. Then, by definition of the
eccentricity of the norm $\|\cdot\|_{\setA}$ with respect to
the norm $\|\cdot\|$,
\begin{equation}\label{E-CH}
\|B-A_{B}\|_{\setA}\le C\cdot\|B-A_{B}\|\le C\cdot H(\setA,\setB),
\end{equation}
where $C=\ecc(\|\cdot\|_{\setA},\|\cdot\|)$.

Consider the obvious inequality
\[
\|Bx\|_{\setA}\le\|A_{B}x\|_{\setA}+\|(B-A_{B})x\|_{\setA},\quad
\forall~ x\in\mathbb{K}^{d}.
\]
Here $\|A_{B}x\|_{\setA}\le\rho(\setA)\|x\|_{\setA}$ because
$\|\cdot\|_{\setA}$ is a Barabanov norm for the matrix set
$\setA$, and $\|(B-A_{B})x\|_{\setA}\le C\cdot
H(\setA,\setB)\|x\|_{\setA}$ by inequality (\ref{E-CH}).
Therefore
\[
\|Bx\|_{\setA}\le\left(\rho(\setA)+C\cdot
H(\setA,\setB)\right)\|x\|_{\setA},\quad \forall~ x\in\mathbb{K}^{d},
\]
and, due to arbitrariness of $B\in\setB$,
\[
\sup_{B\in\setB}\|Bx\|_{\setA}\le\left(\rho(\setA)+C\cdot
H(\setA,\setB)\right)\|x\|_{\setA},\quad \forall~ x\in\mathbb{K}^{d}.
\]
From here by Lemma~\ref{L-rhorel}
\begin{equation}\label{E-rBlerA}
\rho(\setB)\le \rho(\setA)+C\cdot H(\setA,\setB).
\end{equation}

Now, let us prove that
\begin{equation}\label{E-rBgerA}
\rho(\setB)\ge \rho(\setA)-C\cdot H(\setA,\setB).
\end{equation}
By definition of the Hausdorff metric, for any matrix
$A\in\setA$ there is a matrix $B_{A}\in\setB$ such that
$\|B_{A}-A\|\le H(\setA,\setB)$. Then, as before,
\begin{equation}\label{E-CH1}
\|B_{A}-A\|_{\setA}\le C\cdot\|B_{A}-A\|\le C\cdot H(\setA,\setB).
\end{equation}

By evaluating with the help of (\ref{E-CH1}) the second summand
in the next obvious inequality
\[
\|B_{A}x\|_{\setA}\ge
\|Ax\|_{\setA}-\|(B_{A}-A)x\|_{\setA},\quad
\forall~ x\in\mathbb{K}^{d},
\]
we obtain:
\[
\|B_{A}x\|_{\setA}\ge \|Ax\|_{\setA} -C\cdot
H(\setA,\setB)\|x\|_{\setA},\quad \forall~ x\in\mathbb{K}^{d}.
\]
Maximizing now the both sides of this last inequality over all
$A\in\setA$ (which is possible due to arbitrariness of
$A\in\setA$), we get:
\[
\sup_{A\in\setA}\|B_{A}x\|_{\setA}\ge \sup_{A\in\setA}\|Ax\|_{\setA} -C\cdot
H(\setA,\setB)\|x\|_{\setA},\quad \forall~ x\in\mathbb{K}^{d}.
\]
Here the left-hand part of the inequality does not exceed
$\sup_{B\in\setB}\|Bx\|_{\setA}$, while by Barabanov's Theorem
$\sup_{A\in\setA}\|Ax\|_{\setA}\equiv\rho(\setA)\|x\|_{\setA}$.
Hence,
\[
\sup_{B\in\setB}\|Bx\|_{\setA}\ge
\left(\rho(\setA)-C\cdot H(\setA,\setB)\right)\|x\|_{\setA},\quad
\forall~ x\in\mathbb{K}^{d},
\]
from which by Lemma~\ref{L-rhorel} we obtain (\ref{E-rBgerA}).

Inequalities (\ref{E-rBlerA}), (\ref{E-rBgerA}) with
$C=\ecc(\|\cdot\|_{\setA},\|\cdot\|)$ imply (\ref{E-Lipecc})
which finalizes the first step of the proof of
Theorem~\ref{T-main}.

\subsection{Proof of estimate (\ref{E-nubound})}\label{S-step2}

By definition of the eccentricity, the quantity
$\ecc(\|\cdot\|_{\setA},\|\cdot\|)$ is defined as follows
\[
\ecc(\|\cdot\|_{\setA},\|\cdot\|)=
\frac{e^{+}(\|\cdot\|_{\setA},\|\cdot\|)}{e^{-}(\|\cdot\|_{\setA},\|\cdot\|)}.
\]
Here, by the definition (\ref{E-eccentr}) of the quantities
$e^{-}(\cdot)$ and $e^{+}(\cdot)$,
\[
e^{-}(\|\cdot\|_{\setA},\|\cdot\|)=\|x^{-}\|_{\setA},\quad
e^{+}(\|\cdot\|_{\setA},\|\cdot\|)=\|x^{+}\|_{\setA}
\]
for some elements $x^{-}$ and $x^{+}$ satisfying $\|x^{-}\|=1$,
$\|x^{+}\|=1$. Hence
\begin{equation}\label{E-eccxx}
\ecc(\|\cdot\|_{\setA},\|\cdot\|)=\frac{\|x^{+}\|_{\setA}}{\|x^{-}\|_{\setA}}.
\end{equation}

By definition of the measure of irreducibility
$\chi_{p}({\setA})$, for elements $x^{-}$ and $x^{+}$ there are
a natural number $m$, matrices  $\tilde{A}_{i}\in{\setA}_{p}$,
$i=1,2,\dots,m$, and real numbers $\lambda_{i}$,
$i=1,2,\dots,m$, such that
\begin{equation}\label{E-maineq}
\chi_{p}({\setA})x^{+}=\sum_{i=1}^{m}\lambda_{i}\tilde{A}_{i}x^{-},\quad
\sum_{i=1}^{m}|\lambda_{i}|\le 1.
\end{equation}
Here each matrix $\tilde{A}_{i}$ is either the identity matrix
or a product of no more than $p$ factors from $\setA$, that is
$\tilde{A}_{i}=A_{i_{k}}\cdots A_{i_{1}}$ with some $k\le p$
and $A_{i_{j}}\in\setA$. If $\tilde{A}_{i}=I$ then
$\|\tilde{A}_{i}\|_{\setA}=1$. If
$\tilde{A}_{i}=A_{i_{k}}\cdots A_{i_{1}}$ then
$\|\tilde{A}_{i}\|_{\setA}\le\left(\rho(\setA)\right)^{k}$
because, by definition of the Barabanov norm,
$\|\tilde{A}_{i_{j}}\|_{\setA}\le\rho(\setA)$ for any matrix
$A_{i_{j}}\in\setA$. Therefore
\[
\|\tilde{A}_{i}\|_{\setA}\le
\max\{1,\left(\rho(\setA)\right)^{k}\}\le
\max\{1,\left(\rho(\setA)\right)^{p}\},
\]
and (\ref{E-maineq}) implies
\[
\chi_{p}({\setA})\|x^{+}\|_{\setA}\le
\max\left\{1,\left(\rho(\setA)\right)^{p}\right\}\|x^{-}\|_{\setA}.
\]
From here and from (\ref{E-eccxx})
\[
\ecc(\|\cdot\|_{\setA},\|\cdot\|)\le
\frac{\max\{1,(\rho(\setA))^{p}\}}{\chi_{p}({\setA})},
\]
and, since $\rho(\setA)\le\|\setA\|$ by (\ref{E-JSR}), this
last inequality implies the estimate (\ref{E-nubound}), which
finalizes the second step of the proof.

The proof of Theorem~\ref{T-main} is completed.

\bibliographystyle{elsarticle-num}
\bibliography{lipjsr}
\end{document}